\documentclass{article}

\usepackage{amsmath,amssymb}
\usepackage{rotating}
\usepackage{hhline}
\usepackage{xypic}
\usepackage{amssymb}
\usepackage{amsmath,amscd}
\usepackage{pst-node,pstricks,multido,pst-plot,pst-text,pst-3d}%
\usepackage{graphicx}
\usepackage{amsfonts}
\newtheorem{example}{Example}[section]

\newtheorem{theorem}[example]{Theorem}

\newtheorem{corollary}[example]{Corollary}

\newtheorem{proposition}[example]{Proposition}

\newtheorem{lemma}[example]{Lemma}

\def\S{{\mathfrak  S}}
\def\cal#1{{\mathfrak #1}}

\def\<{\langle}
\def\>{\rangle}

\def\C{{\mathbb C}}
\def\Z{{\mathbb Z}}
\def\N{{\mathbb N}}
\def\Y{{\mathbb Y}}

\def\d{\partial}

\def\ashuff#1#2#3{
\kern 1pt \vrule height#1 \overline{\vrule height#3 width 0pt
\hskip#2} \rule{.3pt}{#1}\overline{\vrule height#3 width 0pt
\hskip#2} \rule{.3pt}{#1} \kern 1pt }

\def\Det{{\rm Det}}

\def\X{{\mathbb X}}\def\Y{{\mathbb Y}}

\def\Det{{\rm Det}}
\def\Disc#1#2#3{{\cal D}_{#1}(#2;#3)}

\def\sign{\mbox{sign}}

\begin{document}

\author{Adrien Boussicault and Jean-Gabriel Luque\footnote{Université Paris-Est,
Laboratoire d'Informatique de l'Institut Gaspard Monge
(LABINFO-IGM), UMR CNRS 8049, 5 bd Descartes, 77454 Marne-la-Vallée
Cedex 2, France. Adrien.Boussicault,Jean-Gabriel.Luque@univ-mlv.fr}}
\title
{Staircase Macdonald polynomials and the $q$-Discriminant}
\maketitle
 \begin{abstract}
We prove that  a $q$-deformation $\Disc k\X q$ of the  powers  of
the discriminant
  is equal, up to a normalization, to a specialization of a Macdonald polynomial indexed
by a staircase partition. We  investigate the expansion of $\Disc
k\X q$ on different basis of symmetric functions. In particular, we
show that its expansion on the monomial basis can be explicitly
described in terms of standard tableaux and we generalize a result
of King-Toumazet-Wybourne about the expansion of the
$q$-discriminant on the Schur basis.
 \end{abstract}
\section{Introduction}

Let $\X=\{x_1,\dots,x_n\}$ be an alphabet. The $q$-{\it
discriminant}
\[
\Disc 1\X q:=\prod_{i\neq j}(qx_i-x_j),
\]
is a polynomial encountered in different fields of mathematics.
 In
particular, its specialization at $q=1$ is the  discriminant which
is an example of a  symmetric function invariant under the
transformation $x\rightarrow x+1$ and which has been the subject of
many works in invariant theory (by Cayley, Sylvester and MacMahon).

In condensed matter physics, it plays a crucial role in the context
of the fractional quantum Hall effect. Laughlin \cite{Lau} described
it through a wavefunction whose expression involves an even power of
the Vandermonde determinant

\[
\Psi^k_{\rm Laughlin}(\X)=\pm{\Disc 1\X 1}^k\Psi_{\rm
Laughlin}^0(\X).
\]

 In this paper, we give the links between the
 $q$-discriminant and the Macdonald polynomials. More precisely, our main result is that the
 ``polarized powers'' of the $q$-discriminant

\[
\Disc k\X q:=\prod_{l=1}^k\Disc 1\X {q^{2l-1}},
\]
appear when one evaluates some specialization of ``staircase''
Macdonald polynomials.

 The powers of the discriminant ($q=1$)  are encountered  also in the context of generalizations of the Selberg
 integral \cite{Kan,Kor,Sel}.
 These integrals are closely related to the notion of Hankel
hyperdeterminant \cite{LT1,LT2} and Jack polynomials
\cite{Jack1,Jack2}. The Selberg integral admits  $q$-analogue
involving the $q$-discriminant (see {\it e.g.} \cite{Macdo} ex3
p374). It is interesting to remark that such integrals are related
to Macdonald polynomials \cite{Warnaar}.

More generally, the specializations $t^aq^b=1$ rise deeper
identities related to the generalization of the Izergin and Korepin
determinant due to Gaudin \cite{LascGaud}.

 The paper is organized as follow. In Section \ref{Not}, we recall notations and properties related to symmetric functions.
  Section \ref{Stair} is devoted to the main theorem of the paper. We prove that the polynomial $\Disc k\X q$ is
  a staircase Macdonald polynomial for a specialization of the parameters $q$ and $t$.
   As an application, in Section \ref{Hank}, we give a formula for the coefficients
   arising
  in the expansion of an even power of the Vandermonde determinant in terms of monomial
  functions. Finally, in Section \ref{Schur}, we generalize a
  theorem of King {\it et al.} about the expansion of the
  $q$-discriminant in terms of Schur functions.

\section{Background and notations\label{Not}}
\subsection{Symmetric functions}
We consider the $\C[[q,t,q^{-1},t^{-1}]]$-algebra $Sym$ of symmetric
functions over an alphabet $\X$, {\it i.e.}  the functions which are
invariant under permutations of commuting indeterminates called
letters. There exists various families of such functions. We shall
need the  generating series of complete function:
\[
\sigma_z(\X):=\sum_{i}S^i(\X)z^i=\prod_{x\in\X}\frac1{1-xz}.
\]

This notation is compatible with the sum $\X+\Y$ and the product
$\X\Y:=\sum_{x\in\X, y\in\Y} xy$  in the following sense
 \[\sigma_z(\X+\Y)=\sigma_z(\X)\sigma_z(\Y)=\sum_{i}S^i(\X+\Y)z^i\]
(see {\it e.g.} \cite{Lasc} 1.3 p 5), and
\[
\sigma_t(\X\Y)=\sum_iS^i(\X\Y)t^i=\prod_{x\in\X}\prod_{y\in\Y}{1\over1-xyt}
\]
(see {\it e.g.} \cite{Lasc} 1.5 p13).
 In particular, if $\X=\Y$ one has $\sigma_z(2\X)=\sigma_z(\X)^2$.
 This definition can be extended for any complex number $\alpha$ by
 putting
 $\sigma_z(\alpha\X)=\sigma_z(\X)^\alpha$.

We will use the Schur basis whose elements $S_\lambda$ are indexed by decreasing partitions and defined by
\[
S_\lambda:=\det\left(S^{\lambda_i-i+j}\right)_{1\leq i,j\leq n},
\]
see {\it e.g.} \cite{Macdo} I.3.4 p41 and \cite{Lasc} 1.4.2 p8.
\subsection{Macdonald Polynomials}
The Macdonald polynomials $(P_\lambda(\X;q,t))_\lambda$ form the
unique basis of symmetric functions orthogonal for the standard
$q,t$ deformation of the usual scalar product on symmetric functions
(see {\it e.g.} \cite{Macdo} VI.4 p322), verifying
\begin{equation}
P_\lambda(\X;q,t)=m_\lambda(\X)+\sum_{\mu\leq\lambda}u_{\lambda\mu}m_\mu(\X).
\end{equation}
where $m_\lambda$ is a monomial function in the notation of \cite{Macdo} I.2.1 p8.
Their generating function is (see {\it e.g.} \cite{Macdo} VI.4.13 p324)
\[K_{q,t}(\X,\Y):=\sigma_1\left({1-t\over 1-q}\X\Y\right)=\sum_\lambda P_{\lambda}(\X;q,t)Q_{\lambda}(\Y;q,t),\]
where $Q_\lambda(\X;q,t)= b_\lambda(q,t) P_\lambda(\Y;q,t)$
with
\[
b_{\lambda}(q,t)=\prod_{(i,j)\in\lambda}
{1-q^{\lambda_i-j+1}t^{\lambda'_j-i}\over
1-q^{\lambda_i-j}t^{\lambda'_j-i+1}},
\]
see {\it e.g.} \cite{Macdo} VI.6.19 p339.

Alternatively, when $\X=\{x_1,\dots,x_n\}$ is a finite alphabet, the
Macdonald polynomials can be defined as the eigenfunctions of the
Sekiguchi-Debiard operator ${\cal M}_1$ (see {\it e.g} \cite{Macdo}
VI.3 p315 and VI.4 p325). Indeed,
\begin{equation}\label{M1}
P_\lambda(\X;q,t){\cal M}_1=[[\lambda]]_{q,t}P_\lambda(\X;q,t),
\end{equation}
where, for any $v\in\N^n$, $[[v]]_{q,t}$ is defined as
\begin{equation}\label{eigenvalues}
[[v]]_{q,t}:=q^{v_1}t^{n-1}+q^{v_2}t^{n-2}+\dots+q^{v_n}.
\end{equation}

This operator may be defined in terms of divided differences
\begin{equation}
f(\X){\cal M}_1=f(\X-(1-q)x_1)R(tx_1;\X-x_1)\partial_1\dots
\partial_{n-1}.
\end{equation}
where, for each $i=1\dots n-1$,  $\partial_i$, denoted on the right,
is the operator (see {\it e.g.} \cite{LLM})
\begin{equation*}
f(x_1,\dots,x_n)\partial_i:={f(x_1,\dots,x_i,x_{i+1},\dots,x_n)
-f(x_1,\dots,x_{i+1},x_{i},\dots,x_n)\over
x_i-x_{i+1}}.\end{equation*}

\section{Staircase Macdonald polynomials \label{Stair}}
Let us denote by $\rho:=[n-1,\dots,1,0]$  and set
$m\rho:=[m(n-1),\dots,m,0]$ for $m\in\N$. We  need the following
lemma. 
\begin{lemma}\label{dimeigenspace}
Under the specialization  $t\rightarrow q^{(1-2k)\over 2}$, the
Macdonald polynomial $P_{2k\rho}(\X;q,q^{(1-2k)\over 2})$ belongs to
an eigenspace of ${\cal M}_1$ whose dimension is $1$ and its
associated eigenvalue is
\begin{equation}\label{eigenspec}
[[2k\rho]]_{q,q^{1-2k\over
2}}=\sum_{i=1}^nq^{(2k+1)(n-1)/2}.\end{equation}
\end{lemma}
{\bf Proof} From Equation (\ref{eigenvalues}), the eigenvalue
associated to a  partition $\lambda$ is
 $$ [[\lambda]]_{q,q^{1-2k\over
2}}=\sum_{i=1}^nq^{(1-2k)(n-i)/2+\lambda_i}.$$ Then, if
$[[\lambda]]_{q,q^{1-2k\over 2}}=[[2k\rho]]_{q,q^{1-2k\over 2}}$, it
exists a permutation $\sigma\in\S_n$ such that, for each $1\leq
i\leq n$, one has $
\frac12(2k+1)(n-\sigma(i))=\frac12(1-2k)(n-i)+\lambda_{i}. $ It
follows that
\begin{equation}\label{lambda}\lambda_{i}-\lambda_{i+1}=\frac12(2k+1)(\sigma(i+1)-\sigma(i))-\frac12(1-2k).
\end{equation}
 Since $\lambda$ is a partition, one has necessarily
$\lambda_{i}-\lambda_{i+1}\geq 0$ and Equality (\ref{lambda})
implies $\sigma(i+1)-\sigma(i)\geq{1-2k\over 1+2k}>-1$. This implies
that $\sigma$ is the identity and $\lambda=2k\rho$.
 $\Box$

For simplicity, we set $p:=q^{-\frac12}$ and  we will consider a
finite alphabet $\X=\{x_1,\cdots,x_n\}.$
  Our main result is that
the polarized powers ${\cal D}_k(\X,p)$ of the discriminant are staircase Macdonald
polynomials for the specialization considered here.
\begin{theorem}\label{D2P}
One has
\begin{equation}
\Disc k\X p=(-p)^{\frac12k^2n(n-1)}P_{2k\rho}(\X;q,p^{2k-1}).
\end{equation}
\end{theorem}
{\bf Proof} Reordering factors in $\Disc k{qx_1,x_2,\dots,x_n}p
R(p^{2k-1}x_1;\X-x_1)$, one obtains
\begin{multline}\label{DR2DReq}\displaystyle
\Disc k{qx_1,x_2,\dots,x_n}pR(p^{2k-1}x_1;\X-x_1)=
\displaystyle\Disc k\X pR({p^{-(2k+1)}}x_1;\X-x_1).
\end{multline}
Hence, applying  Equation (\ref{DR2DReq}), the polynomial $\Disc k\X
p{\cal M}_1$ can be rewritten as
\[
\Disc k\X p{\cal M}_1=\Disc k\X pR(p^{-2k-1}x_1;\X-x_1)\d_1\cdots
\d_{n-1}.
\]
Since the polynomial $\Disc k\X p$ is symmetric in $\X$, it commutes
with $\d_1,\dots, \d_{n-1}$ and then
\[
 \Disc k\X p R(p^{-2k-1}x_1;\X-x_1)\d_1\cdots \d_{n-1}= R(p^{-2k-1}x_1;\X-x_1)\d_1\cdots \d_{n-1}\Disc k\X p.
 \]

The remaining factor $R(p^{-2k-1}x_1;\X-x_1)$ is of total degree
$n-1$ and therefore is sent to a constant under $\d_1\dots\d_{n-1}$.
We use the following lemma to compute this constant. 
\begin{lemma}\label{result}
For any letters $a,b$,
\begin{equation} \label{res}
R(ax_1;bx_2,\cdots,bx_n)\d_1\cdots \d_{n-1}=\sum_{i+j=n-1}a^ib^j.
\end{equation}
\end{lemma}
{\bf Proof} Rewrite  $R(ax_1;bx_2,\dots,bx_n)$ as
$$S_{n-1}(ax_1-b(\X-x_1))=S_{n-1}((a+b)x_1-b\X)=\sum x^i_1S_i(a+b)S_{n-1-i}(-b\X).$$
 The image of this sum under $\partial_1\dots\partial_{n-1}$ is
 $S_{n-1}(a+b)S_0(-\X)$ as wanted.
 $\Box$
\\ \\
Applying Lemma \ref{result}, one obtains the value of $\Bbbk$,
\begin{equation}
\Bbbk=\displaystyle
\sum_{i=1}^np^{(2k+1)(i-n)}=\sum_{i=1}^np^{(2k-1)(n-i)-4k(n-i)}.
\end{equation}
From Equality (\ref{eigenspec}), one recognizes that
$
\Bbbk =[|2k\rho|]_{q,p^{2k-1}}.
$
This shows that
\begin{equation}
{\cal D}_k(\X;p)=\beta_{k,n}(p)P_{2k\rho}(\X;q,p^{2k-1}),
\end{equation}
where $\beta_{k,n}(p)$ is a constant  depending only on $p$, $k$ and
$n$. It remains to compute the coefficient $\beta_{k,n}(p)$.
 Since we know that the dominant coefficient
 in
 $P_{2k\rho}(\X;q,p^{2k-1})$ is $1$ by definition,
 it suffices to compute
 the coefficient of the monomial $x_n^{2k(n-1)}\cdots x_2^{2k}$ in ${\cal D}_k(\X,p)$. One finds
\[\beta_{k,n}(p)=(-p)^{\frac12k^2n(n-1)}.\]
This ends the proof.$\Box$
\begin{example}
For $k=2$ and $n=4$, one obtains
\[
P_{[12\ 840]}(x_1+x_2+x_3+x_4;q,q^{3/2})=
{q^{12}}\prod_{i\neq j}\left((qx_i-x_j)(q^3x_i-x_j)\right)
\]
\end{example}
%
\section{Expansion of Macdonald polynomials in terms of  monomial
functions\label{Hank}} Macdonald gives in \cite{Macdo} VI.7.10 p345
the following expansion of the polynomials $Q_\lambda$ in terms of
monomial functions:

\begin{equation}\label{Q2m}
Q_{\lambda}=\sum_\mu\left(\sum_{T}\phi_T(q,t)\right)m_\mu,
\end{equation}
where the inner sum is over the tableaux  of shape $\lambda$ and
evaluation $\mu$ and each $\phi_T(q,t)$ is an explicit  rational
function given in \cite{Macdo} VI.7.11 p346.

 Theorem \ref{D2P} and Equality (\ref{Q2m}) furnish an
expansion of  $\Disc k\X p$ according to the monomial basis,
\begin{equation}\label{D2m}
    {\cal D}_k(\X;p)=
    {(-p)^{\frac12k^2n(n-1)}\over
     b_{2k\rho}(q,p^{2k+1})}\sum_{\lambda}\left(\sum_{T}\phi_T(q,p^{2k-1})\right)m_\lambda
\end{equation}
where the inner sum is over the tableaux of shape $2k\rho$ and
evaluation $\lambda$.

Recall that Jack polynomials \cite{Jack1,Jack2} $P_\lambda^{(\alpha)}(\X)$ are obtained from $P_\lambda(\X;q,t)$ setting $q=t^\alpha$
 and taking the limit when $t$ tends to $1$ (see \cite{Macdo} VI 10).
One has
\begin{equation}
P_\lambda^{(\alpha)}(\X)=\lim_{t\rightarrow1}P_\lambda(\X;t^\alpha,t),
\end{equation}
and
\begin{equation}
Q_\lambda^{(\alpha)}(\X)=\lim_{t\rightarrow1}Q_\lambda(\X;t^\alpha,t)=b_\lambda^{(\alpha)}P_\lambda^{(\alpha)}
\end{equation}
where
$ b_\lambda^{(\alpha)}:=\lim_{t\rightarrow1}b_\lambda(t^\alpha,t).$
Putting
$$\phi_T^{(\alpha)} := \lim_{t\rightarrow 1}\phi_T(t^\alpha,t),
$$
one get from Equation (\ref{D2m}) an expansion of integral powers of
the discriminant.
\begin{corollary}\label{pV2m}
One has
\begin{equation}
\begin{array}{rcl}
\Disc 1\X 1^{k} =\Disc k\X 1  &=&   (-1)^{kn(n-1)\over 2}P_{2k\rho}^{(\alpha_k)}(\X)\\
                  &=&  \displaystyle(-1)^{kn(n-1)\over 2}\left( b^{(\alpha_k)}_{2k\rho}\right)^{-1}
                        \sum_\lambda \left(\sum_T \phi^{(\alpha_k)}_T\right)m_\lambda,
\end{array}
\end{equation}
where $\alpha_k={-2\over 2k-1}$ and the inner sum is over the
tableaux of shape $2k\rho$ and evaluation $\lambda$.
\end{corollary}

\begin{example}\rm
Consider an alphabet $\X=\{x_1,x_2,x_3\}$ of size $3$. One has,
{\footnotesize\begin{eqnarray*} Q_{42}(\X;q,t)=
        \begin{array}{|c|c|c|c|} \hhline{--} 2&2\\\hline 1&1&1&1\\\hline
        \end{array}m_{42}
 +      \begin{array}{|c|c|c|c|} \hhline{--} 2&3\\\hline
        1&1&1&1\\\hline
        \end{array}\ m_{411}
 +
        \begin{array}{|c|c|c|c|} \hhline{--} 2&2\\\hline 1&1&1&2\\\hline
        \end{array}
m_{33}\\
 + \left(
        \begin{array}{|c|c|c|c|} \hhline{--} 2&3\\\hline
        1&1&1&2\\\hline
        \end{array}
 +       \begin{array}{|c|c|c|c|} \hhline{--} 2&2\\\hline
        1&1&1&3\\\hline
        \end{array}
\  \right)m_{321}\\
+ \left(\
        \begin{array}{|c|c|c|c|} \hhline{--} 3&3\\\hline
        1&1&2&2\\\hline
        \end{array}
+       \begin{array}{|c|c|c|c|} \hhline{--} 2&3\\\hline
        1&1&2&3\\\hline
        \end{array}
+       \begin{array}{|c|c|c|c|} \hhline{--} 2&2\\\hline
        1&1&3&3\\\hline
        \end{array}\ \right)m_{222}.
\end{eqnarray*}}
Each tableau $T$ is interpreted as the function $\Phi_T$,
{\footnotesize\begin{eqnarray*}
Q_{42}(\X;q,t)&=&\left(\frac{1-t}{1-q}\right)^2\left(\frac{1-tq}{1-q^2}\right)^2
\left(\frac{1-t^2q^2}{1-tq^3}\right)\left(\frac{1-t^2q^3}{1-tq^4}\right)m_{42}\\
&&+ \left(\frac{1-t}{1-q}\right)^3\left(\frac{1-tq}{1-q^2}\right)\left(\frac{1-t^2q^3}{1-tq^4}\right)
\left(\frac{1-t^2q^2}{1-q^3t}\right)m_{411}+\dots
\end{eqnarray*}
}

Setting $q=t^{-2}$ and taking the limit $t\rightarrow 1$, the
algorithm described here allows to compute the expansion of the Jack
polynomials according to the monomial functions. After
simplification, one obtains
\[ Q_{42}^{(-2)}(\X)={\frac
{1}{280}}\,m_{{4,2}}-{\frac {1}{140}}\,m_{{4,1,1}}-{\frac {1}{
140}}\,m_{{3,3}}+{\frac {1}{140}}\,m_{{3,2,1}}-{\frac
{3}{140}}\,m_{{2,2,2}}.
\]
And finally,
\[
\Disc1\X1=-m_{{4,2}}+2\,m_{{4,1,1}}+2\,m_{{3,3}}-2\,m_{{3,2,1}}+6\,m_{{2,2,2}}.
\]
\end{example}

Corollary \ref{pV2m} can be applied to expand Hankel
hyperdeterminants. Hyperdeterminants are polynomials defined by
Cayley in the aim of generalizing the notion of determinant to
higher dimensional arrays\footnote{Note that Cayley proposed several
generalizations of determinants. The polynomial considered here is
the simplest one in the sense that it generalizes the expansion of
determinant as an alternated sum. Reader can refer to
\cite{LT1,LT2,Matsu,Rice,So1} for more informations on the subject.
} \cite{Ca0,Ca1}. Given a $m$th order tensor ${\rm
M}=\left(M_{i_1\dots i_m}\right)_{1\leq i_1,\dots,i_m\leq n}$ on a
$n$ dimensional space, its hyperdeterminant is
\[
\Det({\rm
M})=\frac1{n!}\sum_{\sigma_1,\dots,\sigma_m\in\S_n}\sign(\sigma_1\dots\sigma_m)\prod_{i=1}^mM_{\sigma_1(i)\dots\sigma_m(i)}.
\]
Note that this polynomial vanishes when $m$ is odd. Suppose that $m=2k$ is an even integer.
 An Hankel hyperdeterminant is an hyperdeterminant whose entries depend only on the sum of the
 indices $M_{i_1\dots i_{2k}}=f(i_1+\dots+i_{2k})$. This kind of hyperdeterminant have been already considered by the authors
  in collaboration with Thibon and Belbachir \cite{LT1,LT2,BBL}.
 In particular, it is shown that the coefficients $C_\lambda(n,l)$ arising in the expression
\[
\Det\left(M_{i_1+\dots+i_{2k}}\right)=\sum_\lambda C_\lambda(n,k)\prod_{i=1}^nf(\lambda_i),
\]
are equal (up to  a multiplicative term equal to the number of
permutations of $\lambda$ divided by $n!$) to those arising in the
expansion of $\Disc k\X 1$ in terms of monomial functions.
\begin{example}\rm
From the expansion of the Jack polynomial $P_{84}^{(-2/3)}$, for an alphabet of size $3$,
\[
\begin{array}{rcl}
P_{84}^{(-2/3)}(x_1+x_2+x_3)&=& m_{84}-4m_{831}+6m_{822}-4m_{75}+12m_{741}-8m_{732}
+6m_{66}-8m_{651}\\&&-22m_{642}+48m_{633}+48m_{552}-36m_{543}+90m_{44},
\end{array}
\]
one deduces the expansion of the Hankel hyperdeterminant

\[\begin{array}{rcl}
\Det\left(f(i_1+i_2+i_3+i_4)\right)_{0\leq i_1,i_2,i_3,i_4\leq 3}&=& f(8)f(4)f(0)-4f(8)f(3)f(1)+3f(8)f(2)^2\\
&&- 4f(7)f(5)f(0)+12f(7)f(4)f(1)-8f(7)f(3)f(2)
\\&&+3f(6)^2f(0)-8f(6)f(5)f(1)-22f(6)f(4)f(2)\\&& +24f(6)f(3)^2+24f(5)^2f(2)-36f(5)f(4)f(3)\\&&+15f(4)^3.
\end{array}
\]

\end{example}
Furthermore, in \cite{LLM} Lapointe {\it et al.} gave a
determinantal expression of Jack polynomial in terms of monomial
functions. These computations leads naturally to a determinantal
expression for Hankel hyperdeterminants.

Note that the formula for the Macdonald polynomials $\tilde
H_\lambda$, given by Haglund, Haiman and Loehr \cite{HHL}, provides
an expansion of $\Disc k\X q$ in terms of modified monomial
functions $m_\lambda(\X(1-t))$ having a combinatorial
interpretation.
\section{Expansion of the polarized powers of the $q$-discriminant in terms of Schur functions\label{Schur}}
Di Francesco {\it et al.} \cite{DGIL} considered the problem of the
expansion of the discriminant in terms of Schur functions. They
defined the $n$-{\it admissible} partitions to be the partitions in
the interval $[(n-1)^n]$, $[2(n-1),\dots,2,0]$ (with respect to the
dominance order). They conjectured that they are exactly those
 occurring in the expansion of the
discriminant. This conjecture is false as shown by Scharf {\it et
al.} \cite{STW}. However, Kind {\it et al.} \cite{KTW} proved that
it becomes true when replacing the discriminant by the
$q$-discriminant.

In this section, we generalize this property to ${\cal D}_k(\X;q)$.
We define $(n,m)$-{\it admissible} partitions to be the partitions
which appear in the expansion
\begin{equation}\label{DefAltAdm}
m_\rho(\X)^{m-1}S_\rho(\X)=\sum_\lambda b_\lambda^{n,m}m_\lambda(\X)
\end{equation}
where $\X$ is an alphabet of size $n$. When $m=2k$ is even, the
$(n,2k)$-admissible partitions are those of the interval
$[(k(n-1))^n]$, $[2k(n-1),\dots,2k,0]$. We prove that a partition
appear in the expansion of ${\cal D}_k(\X;q)$ in terms of Schur
functions if and only if it is a $(n,2k)$-partition.

\subsection{Computing admissible partitions}
Let us denote by $A_{n,m}$ the set defined recursively by
\begin{equation}\label{RecA}\begin{array}{l}
A_{n,1}:=\{\lambda=[\lambda_1,\dots,\lambda_n]| \rho\geq\lambda\}\\
A_{n,m}:=\{((\lambda_1+\sigma(1)-1,\dots,\lambda_n+\sigma(n)-1))|\sigma\in\S_n\mbox{
and } \lambda   \in A_{n,m-1}\}.\end{array}
\end{equation}

\begin{lemma}\label{AareAdm} Let $\lambda$ be a partition.
The following assertions are equivalent.
\begin{enumerate}
\item The partition $\lambda$ belongs to $A_{n,m}$.
\item The partition $\lambda$ is $(n,m)$-admissible.
\item $\lambda$ is partition of length $n$ less or equal to $m\rho$ with respect to the
dominance order.
\end{enumerate}
\end{lemma}
{\bf Proof} The equivalence between the assertions 1 and 2  is
straightforward from Equations (\ref{DefAltAdm}) and (\ref{RecA}).
Furthermore, from Equation (\ref{RecA}), the maximal partition of
$A_{n,m}$ is $m\rho$. It remains to prove $3\Rightarrow 1$. We
proceed by induction on $m$, if $m=1$ then the result is trivial.
Suppose that $m>1$. Let $\lambda$ be a partition of size $n$ less or
equal to $m\rho$ with respect to the dominance order. Then
$((\lambda-\rho))$ is a partition less or equal to $(m-1)\rho$.
Indeed, putting $((\lambda-\rho))=(\mu_1,\dots,\mu_n)$, for a
permutation $\sigma\in\S_n$, one has
$\mu_i=\lambda_{\sigma(i)}+n-\sigma(i)$. Hence, for each $i$
\[
\begin{array}{rcl}
\mu_1+\dots+\mu_i&\leq&\lambda_{\sigma(1)}+\dots+\lambda_{\sigma(i)}+n-\sigma(1)+\dots
+n-\sigma(i)\\
&\leq&\lambda_1+\dots+\lambda_i+n-1+\dots+n-i\\
&\leq&(m-1)(n-\frac12i(i+1))
\end{array}
\]
implies  $((\lambda-\rho))\leq (m-1)\rho$ for the dominance order.

 By
induction, $((\lambda-\rho))$ belongs to $A_{n,m-1}$. Furthermore,
it exists a permutation $\sigma$ such that
$((\lambda-\rho))+\rho^\sigma=\lambda$. Hence, from Equation (\ref{RecA}),
$\lambda\in A_{n,m}$.$\Box$


\subsection{Counting  admissible partitions}

One considers the free commutative monoid ${\cal T}$ generated by
the symbols $T=\{\tau_1,\dots, \tau_{n-1}\}$ acting on the vectors
of size $n$ by
\[
\tau_i[v_1,\dots,v_n]=[v_1,\dots,v_{i-1},v_i-1,v_{i+1}+1,v_{i+1},\dots,v_n].\]
For a given vector $v\in\Z^n$, ${\cal T}.v$ is the set of the
vectors $w=[w_1,\dots,w_n]\in\Z^n$ of same weight ({\it i.e.}
$v_1+\dots+v_n=w_1+\dots w_n$) lower or equal to $v$ for the
dominance order. In particular, if $v=\lambda$ is a partition then
${\cal T}.\lambda$ contains all the partition of size $n$ lower or
equal to $\lambda$. To each vector $v\in\Z^n$, one associates the
monomial $z^{v}=z_1^{v_1-v_2}\dots z_{n-1}^{v_{n-1}-v_n}$. For a
given weight, the monomial $z^{v}$ characterizes completely $v$,
furthermore $v$ is a (decreasing) partition if and only if its
weight is non negative and the degree of the monomial $z^v$ in each
variable $z_i$ is non-negative.
\begin{example}
\begin{equation}
\begin{array}{cccccccc}
         &        &                        & {z_1^3\atop [4,1,1]}     &                      &   \\
         &        &   \tau_1\swarrow       &             &   \searrow\tau_2     &\\
         &        &         {z_1z_2\atop[3,2,1]}        &             & {z_1^4,z_2^{-2}\atop [4, 0,2]}            & \\
        &\swarrow&                         &\searrow\swarrow &                      &\searrow &\\
         &{ z_1^{-1}z_2^2\atop[2,3,1]} &
&{z_1^2z_2^{-1}\atop[3,1,2]}  &
                               &{z_1^4z_2^{-3}\atop [3,-1,2]}\\
\swarrow &         &         \searrow \swarrow     &         &\searrow\swarrow&  &\searrow \\
\ {z_1^{-3}z_2^3\atop[1,4,1]}&        &      {1  \atop [2,2,2]} &
& {z_1^3z_2^{-3}  \atop     [3,0,3]}       &            &{x_1^{-5}z_2^5\atop [3,-2,3]}\\
         &        &                        &\dots    &                      &            &
        \end{array}
\end{equation}
\end{example}
Acting on $v$ by $\tau_i$ is equivalent to multiply $z^v$ by
\[
t_i=\left\{\begin{array}{ll}{z_{i-1}z_{i+1}\over z_i^2}&\mbox{ if
} 1<i<n-1\\
{z_{2}\over z_1^2}&\mbox{ if } i=1\\
{z_{n-2}\over z_{n-1}^2}&\mbox{ if } i=n-1.
\end{array}\right.
\]
Since there is no algebraic relations between the $t_i's$, each
vector appears in ${\cal T}.v$ with  multiplicity $0$ or $1$. In
other words, one has
\begin{equation}\label{EqSeries}
\sigma_q(T).z^v=\prod_i\frac1{1-tq}.z^v=\sum_{w\leq v}
q^{\alpha_{v,w}}z^w.
\end{equation}
where ${\alpha_{v,w}}$ is the degree of the monomial acting on $v$
to obtain $w$. Extracting the monomial which encodes a partition is
equivalent to extract the part of the series (\ref{EqSeries})
constituted only with non-negative exponents. This operation is
performed by the MacMahon Omega operator  (see {\it e.g}
\cite{Andrews})
\[
\Omega_{x_1,\dots,x_p} \sum_{n_1,\dots,n_p\in\Z}
\alpha_{n_1,\dots,n_p}x_1^{n_1}.\dots.x_p^{n_p}=\sum_{n_1,\dots,n_p\in\N}
\alpha_{n_1,\dots,n_p}x_1^{n_1}.\dots.x_p^{n_p}.
\]

\begin{example}\rm
 One has
\[\Omega_{z_1,z_2}\frac{z_1^3}{(1-\frac{z_1}{z_2^2}q)(1-\frac{z_2}{z_1^2}q)}=z_1^3+qz_1z_2+q^3 \]
which implies that the set of the partitions of size $3$ lower or
equal to $[411]$ is $\{[411],[321],[222]\}$.
\end{example}
Hence,
\begin{proposition}
The size $n\geq 2$ of the alphabet being fixed, the generating
series of the $(n,k)$-admissible partitions is the rational function
\[
{\cal
A}_{n}(q,t;z_1,\dots,z_{n-1})=\Omega_{z_1,\dots,z_{n-1}}\left((1-tz_1\dots
z_{n-1})(1-\frac{z_0z_2}{z_1^2}q)(1-\frac{z_1z_3}{z_2^2}q)\dots
(1-\frac{z_{n-2}z_n}{z_{n-1}^2}q)\right)^{-1},
\]
where $z_0=z_n=1$.
\end{proposition}
\begin{example} \rm
Let us give the first value of ${\cal
A}_{n}(q,t;z_1,\dots,z_{n-1})$.
\begin{enumerate}
\item First, one considers the special case $n=2$,
$$\begin{array}{rcl}
{\cal A}_{1}(q,t;z_1)&=&\Omega_{z_1}\left((1-tz_1)(1-\frac
q{z_1^2})\right)^{-1}=\left((1-qt^2)\left (1-{\it z_1}\,t
\right)\right)^{-1}\\&=&1+z_1t+(q+z_1^2)t^2+(qz_1+z_1^3)t^3+(q^2+qz_1^2+z_1^4)t^4+\dots.
\end{array}
$$
This means that, for $k=1$ the only admissible partition is $[21]$,
for $k=2$, there is two admissibles partitions $[42]$ and $[33]$,
for $k=3$ the admissibles partitions are $[63]$ and $[54]$ etc...
\item If $n=3$,
$ {\cal A}_{3}(q,t;z_1,z_2)=-{\frac {1-{{\it z_1}}^{3}{{\it
z_2}}^{3}{q}^{2}{t}^{4}}{\left (1-t{ \it z_1}\,{\it z_2}\right
)\left (1-q{{\it z_2}}^{3}{t}^{2}\right )\left (1-{{\it
z_1}}^{3}q{t}^{2}\right )\left (1-{q}^{2}t\right )}}.
$
\item If $n=4$,
$ {\cal A}_{4}(1,t;1,1)={\frac
{{t}^{4}+5\,{t}^{3}+7\,{t}^{2}+2\,t+1}{\left (1+t\right )^{2} \left
(1-t\right
)^{4}}}
.
 $
 \item If $n=5$,
$
{\cal A}_{5}(1,t;1,1)={\frac
{3\,{t}^{6}+21\,{t}^{5}+61\,{t}^{4}+68\,{t}^{3}+39\,{t}^{2}+7\,t+1}{\left
(1-t\right )^{5}\left (t+1\right
)^{3}}}
$
\end{enumerate}
\end{example}

\subsection{Characterization of the partitions arising in the expansion of $\Disc k\X q$}

In this paragraph, one extends the result of King-Toumazet-Wybourne
to the polynomials ${\cal D}_k(\X;q)$.

\begin{theorem}\label{Adm}
Expand $\Disc k\X q$ in terms of Schur functions, $$\Disc k\X
q=\sum_\lambda c_\lambda(q)S_\lambda(\X).$$ Then, $c_\lambda(q)\neq
0$  if and only if $\lambda$ is a $(n,2k)$-admissible partition.
\end{theorem}
{\bf Proof} Let us prove first the only if part. From Theorem
\ref{D2P}, the polynomial $\Disc k\X q$  equals (up to a
multiplicative coefficient)  a specialization of the Macdonald
$P_{2k\rho}(\X;q,t)$. But it is well known that the partitions
arising in the expansion of $P_{2k\rho}(\X;q,t)$ in terms of Schur
functions belong to the interval $[(k(n-1))^n]$, $2k\rho$ (see {\it
e.g.} the determinantal expression of Macdonald polynomials given in
\cite{LLM}). From  Lemma \ref{AareAdm}, this is equivalent to the
fact that $\lambda$ is $(n,2k)$-admissible.

Conversely, to prove that the admissibility of $\lambda$ implies the
non nullity of $c_\lambda(q)$, it suffices to prove it for a
specialization. We will set $q=-1$. In this case,
$$\Disc k\X q=\prod_{i\neq j}(x_i+x_j)^k=S_\rho(\X)^{2k}.$$

We will prove a stronger result  showing that the coefficient
$c_\lambda^{n,m}$ in the expansion
$$S_\rho(\X)^{m}=\sum_\lambda c_\lambda^{n,m}S_\lambda(\X) $$
is non-zero if and only if $\lambda$ is ($n,m)$-admissible. We
proceed  by induction on $m$. Note that the initial case ($m=2$)
have been proved by King-Toumazet-Wybourne in \cite{KTW} Corollary
3.2 as a consequence of an important result of
Bereinstein-Zelevinsky \cite{BZ}.

One needs the two following lemmas
 \begin{lemma}\label{nk2n-1k}
If $\lambda$ is a $(n,m)$-admissible partition ($m>1$), then
$((\lambda-\rho))$ is a $(n,m-1)$-admissible partition.
 \end{lemma}
 {\bf Proof}
From Equality (\ref{DefAltAdm}), each $(n,m)$-admissible partition
can be obtained by adding a permutation of $\rho$ to a
$(n,m-1)$-admissible partition. This is equivalent to our statement.
 $\Box$

\begin{lemma}\label{specLR}
Let $\mu\subset\lambda$ be a partition and
$\nu:=((\lambda_1-\mu_1,\dots,\lambda_{n-1}-\mu_{n-1},\lambda_n-\mu_n))$.
Then, the Littlewood-Richardson coefficient
$c_{\mu\nu}^\lambda=\langle S_\lambda,S_\mu S_\nu\rangle$ equals
$1$.
\end{lemma}
{\bf Proof} The Littlewood-Richardson coefficient $c_{\mu\nu}^\lambda$
 is equal to the number of tableaux of shape $\nu$ and evaluation
 $\lambda-\mu$. But $\lambda-\mu$ is a permutation if $\nu$ and Theorem 11.4.3 of
 \cite{Lasc} implies that such a tableau exists and is unique. This
 ends the proof.
  $\Box$\\ \\
{\bf End of the proof of Theorem \ref{Adm}} Let $\lambda$ be a
$(n,m)$-admissible partition. Since $\rho\subset\lambda$,  Lemma
\ref{nk2n-1k} implies that the partition $\mu=((\lambda-\rho))$ is
$(n,m-1)$-admissible. And by induction, $S_\mu$ appears with a
non-zero coefficient in $S_\rho^{m-1}$. The positivity of the
Littlewood Richardson coefficients implies that each partition $\nu$
such that $c^\nu_{\mu,\rho}\neq 0$ appears with a non-zero
coefficient in the expansion of $S_\rho^{m}$. In particular, from
Lemma \ref{specLR}, it is the case of $\lambda$. This shows that
$c_\lambda^{n,m}\neq 0$ if and only if $\lambda$ is
$(n,m)$-admissible and proves the Theorem.$\Box$.

\noindent Note that other expansion of Macdonald functions can be
found in literature (for example Hall-Littlewood polynomials can be
expanded in terms of plane partitions \cite{Add1}), it should be
interesting to investigate the properties $\Disc k\X q$ which can be
deduced from these expansions.

\noindent{\bf Acknowledgments}

We are grateful to Alain Lascoux for  fruitful discussions.

 \end{document}